\documentclass[12pt]{article}

\usepackage{amssymb,amsmath,amsthm}

\newtheorem{theorem}{Theorem}[section]
\newtheorem{lemma}[theorem]{Lemma}

\newtheorem{corollary}[theorem]{Corollary}

\theoremstyle{definition}

\theoremstyle{remark}

\numberwithin{equation}{section}

\newcommand{\wedgeco}%
{\displaystyle\operatornamewithlimits{\wedge}_{\raise3pt\hbox{,}}}

\newcommand{\ostar}{\odot\kern-6.4pt\ast}

\def\zs#1{_{\lower2pt\hbox{$\scriptstyle#1$}}}

\newcommand{\hb}{\hbar}
\newcommand{\ol}{\overline}
\newcommand{\la}{\lambda}
\newcommand{\pa}{\partial}

\newcommand{\ds}{\displaystyle}

\newcommand{\wh}{\widehat}

\newcommand{\od}{\overset{\text{\rm def}}{=}}

\newcommand{\const}{\operatorname{const}}
\newcommand{\Exp}{\operatorname{Exp}}
\newcommand{\tr}{\operatorname{tr}}
\newcommand{\id}{\operatorname{id}}

\newcommand{\bI}{\mathbf{I}} 
\newcommand{\bP}{\mathbf{P}}
\newcommand{\bS}{\mathbf{S}}

\newcommand{\bR}{\mathbb{R}}
\newcommand{\bZ}{\mathbb{Z}}

\newcommand{\cA}{\mathcal{A}} 
\newcommand{\cC}{\mathcal{C}} 
\newcommand{\cD}{\mathcal{D}} 
\newcommand{\cE}{\mathcal{E}} 
\newcommand{\cF}{\mathcal{F}} 
\newcommand{\cH}{\mathcal{H}} 
 
\newcommand{\cL}{\mathcal{L}} 
\newcommand{\cX}{\mathcal{X}}

\begin{document}

\title{Quantization and Intrinsic Dynamics} 

\author{Mikhail Karasev\thanks{This research was supported 
in part by the Russian Foundation for Basic Research, 
Grant~02-01-00952.}\\
\\
\small Applied Math. Department\\
\small Moscow Institute of Electronics and Mathematics\\
\small Moscow 109028, Russia\\
\small karasev@miem.edu.ru}

\date{}

\maketitle


\begin{abstract}
A dynamical scheme of quantization of symplectic
manifolds is described.
It is based on intrinsic Schr\"odinger and Heisenberg type
nonlinear 
evolutionary equations with multidimensional time running over
the manifold.  
This is the restricted version of the original article 
to be published in ``Asymptotic Methods for Wave and 
Quantum Problems'' (M.~Karasev, ed.), 
Advances in Math. Sci., AMS.


\end{abstract}

\maketitle

\section{Introduction}

The paper deals with constructing irreducible
quantum manifolds, that is, quantizing symplectic manifolds
\cite{1}--\cite{46}.  

Reviewing results of previous works~\cite{12,13}, 
we investigate an {\it intrinsic evolutionary 
differential equation\/} for the
 integral kernel of a quantum
associative product of functions over the symplectic manifold.
The ``time'' variable in this evolutionary equation is
multidimensional and runs over the same symplectic manifold.
We describe the solution of this equation in the semiclassical
approximation, as well as  in the sense of weak asymptotics, by
purely geometric means. 

It is an open question how to solve this equation exactly (and
so to construct exact quantization) for a general symplectic
manifold. However, in some nontrivial classes of examples, 
this occurs to be possible.

The quantization which we deal with is formal in the sense that 
we do not consider analytic conditions on function spaces,
Hilbert norms, etc. Nevertheless, this approach is not
completely formal, it is adapted to the strict quantization
scheme~\cite{14} and directly related to operator
representations. This is a program for the future:
exact operator realization of the this approach over general
symplectic manifolds. The important case, where the
operator representation is possible,
are K\"ahlerian manifolds, in particular, symplectic leaves
in some Poisson manifolds with partial complex structure.
The equations for the $*$-product kernel in this case were
used in~\cite{15,16}.
About examples for the non-K\"ahlerian situation see,
for instance,~\cite{17,17a,18}.

Besides the algebraic and analytic machinery, which we explain
below, there is an important geometric framework around such a
dynamical scheme of quantization. 
In the present paper we would like to pay attention to the fact 
that the old Ether idea, 
or more precisely, the idea of intrinsic dynamics
on the phase space level,
is mathematically fruitful and, moreover, 
automatically arises from the quantization problem. 
Hopefully, on this way, one can realize the early Weyl's 
objectives and his intuitive program of exploring the
continuum ``relationship between a part and the whole'' 
(see \cite{19}).

\section{Basic wave equations}
 
Let us look first at the stationary Schr\"odinger equation for
wave functions $\psi$ and for the energy levels $\la$:
\begin{equation}
\big(\wh{p}^{\,2}+V(\wh{q}\,)\big)\psi=\la \psi,
\tag{2.1}
\end{equation}
where $\wh{q}=q$, $\wh{p}=-i\hb\pa /\pa q$.
Solving this equation is difficult, and the properties of the
wave functions are complicated.
Say, in the semiclassical approximation we have 
\begin{equation}
\psi\sim c_\hb e^{\frac{i}{\hb}S}(\rho+O(\hb)),\qquad c_\hb=\const.
\tag{2.2}
\end{equation}
But, this expression is very far from the asymptotics of
the wave function in the presence of focal or turning points.

On the other hand, the weak limit $\psi^0$ of the function $\psi$ as
$\hb\to0$ is very simple. 
Indeed, if $\la^0$ denotes the limit of $\la$, then Eq.~(2.1)
at $\hb=0$ becomes
$$
V(q)\psi^0(q)=\la^0\psi^0(q).
$$
The solutions are obvious:
\begin{equation}
\psi^0(q)=\delta_{q^0}(q),\qquad \la^0=V(q^0),
\tag{2.3}
\end{equation}
where $q^0$ is arbitrary.
So, {\it the quantum wave functions in the weak limit $\hb=0$
are just $\delta$-functions\/} characterizing the singular
distribution of the probability amplitude of a classical
particle to be localized at a given 
point.

Now let us look at the algebraic picture behind 
Eq.~(2.1), namely, at the Heisenberg commutation relations
\begin{equation}
[\wh{q},\wh{p}\,]=i\hb\bI.
\tag{2.4}
\end{equation}
These relations introduce a noncommutative
product~$*$ into the function algebra over the $q,p$-phase space
as follows: 
\begin{equation}
\wh{f}\od f(\wh{q},\wh{p}\,),\qquad \wh{f}\wh{g}=\wh{f*g}.
\tag{2.5}
\end{equation}

Here we use the Weyl-symmetrized functions in the noncommuting
operators $\wh{q}, \wh{p}$ so that, 
if $f(q,p)=\exp\{i(\alpha q+\beta p)\}$, 
then $\wh{f}=\exp\{i(\alpha \wh{q}+\beta \wh{p}\,)\}$.

If $\hb=0$, then the algebra (2.4) becomes commutative 
and the product $f*g$ becomes the usual commutative product~$fg$
of functions.  
One can represent it via the integral kernel as 
\begin{equation}
(fg)(z)=\iint \delta_x(z)\delta_y(z)f(x)g(y)\,dm(x) dm(y).
\tag{2.6}
\end{equation}
Here the arguments $x,y,z$ run over the $q,p$-phase space 
and $\delta$ denotes the Dirac delta function with respect to a
certain measure~$dm$ on the phase space.

If $\hb\ne0$, then the noncommutative product $f*g$ in
(2.5) can be represented in integral form as well:
\begin{equation}
(f*g)(z)=\frac1{(2\pi\hb)^{2n}}\iint K_{x,y}(z)f(x)g(y)\,dm(x) dm(y),
\tag{2.7}
\end{equation}
where $2n$ is the dimension of the phase space.
The integral kernel here is the noncommutative product of
$\delta$-functions: 
$$
\frac1{(2\pi\hb)^{2n}}K_{x,y}=\delta_x*\delta_y.
$$

From~(2.6) we see that in the weak limit as $\hb\to0$
the integral kernel is equal to the commutative product of delta
functions:  
$$
\lim_{\hb\to0}\frac1{(2\pi \hb)^{2n}}K_{x,y}=\delta_x\delta_y.
$$
This is the same picture as for the wave function $\psi$ in~(2.1)
and its weak limit $\psi^0$~(2.3).
Looking at this analogy, we conclude that 
{\it a natural interpretation of the kernel of the
$*$-product could be a ``wave function'' of something}.

Of what?

In (2.4), (2.5), and (2.7), there are no
particles, no {\it a priori\/} Hamiltonians. Nevertheless, one
can introduce a Hamiltonian of the product $*$ itself by
mimicking the Schr\"odinger evolutionary equation.

Let us first note that the operators $\wh{f}$ in (2.5)
depend on~$f$ linearly and so can be written as integrals
\begin{equation} 
\wh{f}=\frac1{(2\pi\hb)^n}\int f \bS\,dm
\tag{2.8}
\end{equation}
over the phase space. Here $\bS=\{\bS_x\}$ is a family of operators
parametrized by points~$x$ of the phase space. 
In the case of the flat phase space $\bR^{2n}$ one could choose
the Liouville measure, $dm(x)=dx=dq dp$;
then the explicit formula for the integral kernel of the
operators $\bS_x$ in the Hilbert space $L^2(\bR^n)$ is 
$$
\bS_x\sim\delta\bigg(q-\frac{q'+q''}{2}\bigg)
\exp\bigg\{\frac{i}{\hb}p(q'-q'')\bigg\},\qquad x=(q,p).
$$

An important fact about the representation~(2.8) is that it is
easily inverted~\cite{21}:
\begin{equation}
f=\tr(\wh{f} \bS).
\tag{2.8a}
\end{equation}

The integral kernel $K_{x,y}$ in (2.7) is related to $\bS_x$ 
by the formula
\begin{equation}
\wh{K}_{x,y}=\bS_x\bS_y.
\tag{2.9}
\end{equation}

It is not surprising that each $\bS_x$ is self-adjoint
(if one wants to associate self-adjoint operators $\wh{f}$ with
real symbols $f$ by (2.8)), but it is remarkable that all
$\bS_x$ are almost unitary; namely, $\bS_x$ differs from a
unitary only by a constant multiplier
$$
\mu=2^n.
$$
Thus we have 
\begin{equation}
\bS^*=\bS,\qquad \bS^2=\mu^2\cdot\bI.
\tag{2.10}
\end{equation}

The first consequence from this almost unitarity is the
representation 
\begin{equation}
\bS=\mu(\bI-2\bP),
\tag{2.11}
\end{equation}
where $\bP=\{\bP_x\}$ is a family of orthogonal projections.

Second, if one looks at the parameter~$x$ of the almost unitary
family $\bS=\{\bS_x\}$ as at a ``time'' variable, then one can
derive the generator of $\bS$ with respect to this
multidimensional ``time'': 
\begin{equation}
i\hb\pa\bS=\wh{\cH}^\hb\bS.
\tag{2.12}
\end{equation}
Here the differential $\pa=\pa/\pa x$ is taken 
with respect to the variable~$x$ running  over the phase space.
The operator $\wh{\cH}^\hb$ is the self-adjoint generator given
by 
\begin{equation}
\wh{\cH}^\hb=2i\hb(\pa\bP\cdot\bP-\bP\cdot \pa\bP).
\tag{2.13}
\end{equation}

Since $\bS$ and $\wh{\cH}^\hb$ are self-adjoint, 
we conclude from~(2.12) that
\begin{equation}
\wh{\cH}^\hb\bS=-\bS\wh{\cH}^\hb
\tag{2.14}
\end{equation}
and so the equation~(2.12) can be transformed as 
\begin{equation}
i\hb\pa \bS=\frac12[\wh{\cH}^\hb,\bS].
\tag{2.12a}
\end{equation}
This dynamical equation looks like the Heisenberg evolution
equation,  
or like a Lax type equation for the ``$L-A$ pair.''

By means of (2.8a) and (2.9), we return to the integral kernel
$K_{x,y}$ and obtain from (2.12): 
\begin{equation}
i\hb\pa_x K_{x,y}=\cH^\hb_x* K_{x,y}.
\tag{2.15}
\end{equation}
The Cauchy data are
\begin{equation}
K_{y,y}=\mu^2,\qquad \forall\,y.
\tag{2.16}
\end{equation}

The constant $\mu^2$ is determined by the condition
that the unity function~$1$ is 
the unity element for product~(2.7), i.e., the function 
\begin{equation}
\frac1{(2\pi\hb)^{2n}}\int K_{x,y}(z)\,dm(y)\sim \delta_x(z)
\tag{2.17}
\end{equation}
in arguments $x,z$ serves as the kernel of the
unity operator  in the quantum function algebra, that is,
it equals $\delta_x(z)$ for the Euclidean phase
space.

Also from the self-adjointness~(2.10), 
we have 
\begin{equation}
\overline{K_{x,y}(z)}=K_{y,x}(z),\quad\text{or}\quad
\overline{f*g}=\overline{g}*\overline{f}.
\tag{2.18}
\end{equation}

In addition to these simple conditions,
there is also the  fundamental {\it cyclicity condition},
which follows from the trace representation 
$K_{x,y}(z)=\tr(\bS_x\bS_y\bS_z)$, namely,
\begin{equation}
K_{x,y}(z)=K_{z,x}(y),\quad\text{or}\quad
\int f*g\, dm=\int fg\,dm.
\tag{2.19}
\end{equation}
Actually, this is the condition for the choice of the measure~$dm$.

Note that, after integrating (2.15) with the function 
$1(y)f(z)$ and using the property (2.19), one obtains
the identities
\begin{equation}
i\hb df(x)=(f*\cH^\hb_x)(x)=-(\cH^\hb_x*f)(x), 
\qquad \forall\,f.
\tag{2.20}
\end{equation}

These identities determine the ``germ'' of the $*$-product near
a given point~$x$ of the phase space. 
Global information about the $*$-product is contained 
in the wave equation~(2.15), which can be resolved using the
operator multiplicative exponential or the 
$*$-exponential\footnote{Such
exponentials were introduced in early works on deformation
quantization \cite{39,22} and were intensively exploited, for
instance, see~\cite{23} and the references therein.} 
as follows:
$$
\wh{K}_{x,y}=\mu^2\cdot\Exp\bigg\{-\frac{2i}\hb\int^x_y\wh{\cH}^\hb\bigg\},
$$
or
\begin{equation}
f*g=\frac{\mu^2}{(2\pi\hb)^{2n}}\iint f(x) \Exp_*
\bigg\{-\frac{2i}\hb\int^x_y {\cH}^\hb\bigg\}
g(y)\,dm(x) dm(y),
\tag{2.21}
\end{equation}

This formula manifests the dynamical character of quantization.
We see a sum of parallel-transported amplitudes multiplied by
values of functions~$g$ and~$f$ at the initial and final ``time''
points. The role of ``time'' is played by points of the phase
space.

We stress that (2.15), as well as (2.12), are Schr\"odinger type
dynamical equations behind any particles. We see no particles
but, nevertheless, something is ``moving'' and controlled by the
Hamiltonian~$\cH^\hb$.  
We call $\cH^\hb$ the {\it Ether Hamiltonian}.  

Of course, in the case of the standard phase space $\bR^{2n}$,
the kernel~$K$ is well known independently of the calculations
performed above, namely, 
\begin{equation}
K_{x,y}(z)=\exp\bigg\{\frac{i}{\hb}\Phi_{x,y}(z)\bigg\},
\tag{2.22}
\end{equation}
where $\Phi_{x,y}(z)$ is the symplectic area of the triangle
with mid-points $x,y,z$ (see~\cite{24,25} and~\cite{39,26} 
for general flat symplectic spaces).

The Ether Hamiltonian $\cH^\hb$ is easily derived in the case
$\bR^{2n}$ from Eq.~(2.12) or from (2.15), (2.16):
\begin{equation}
\cH^\hb_x(z)=\frac{i\hb}{\mu^2}\pa_x K_{x,y}(z)\bigg|_{y=z}
=2\omega(z-x)\cdot dx.
\tag{2.23}
\end{equation}
Here $\omega$ is the constant matrix of the symplectic form on
$\bR^{2n}$ in the Cartesian coordinates. In the Darboux
coordinates $\omega=J=\left(\begin{matrix}0&I\\-I&0\end{matrix}\right)$.

Even for the simplest integral kernel (2.22) the differential
equation~(2.15) seems to be new. Other attempts to find some
differential equations for the $*$-product kernel were initiated
in \cite{27,28}, but they are certainly different from our approach.
\medskip

Note that the operator $\wh{\cH}^\hb$ in (2.12) is defined via
(2.8), and so Eqs.~(2.12), (2.12a) are actually nonlinear:
\begin{equation}
i\hb\pa \bS_x=\int \cH^\hb_x(y)\bS_y\bS_x\,dy,
\tag{2.24}
\end{equation}
or
\begin{equation}
i\hb\pa \bS_x=\frac12\int \cH^\hb_x(y)[\bS_x,\bS_y]\,dy.
\tag{2.24a}
\end{equation}

Similarly, the product $*$ in the evolutionary equation (2.15)
for the ``wave function''~$K$ is determined by the same
function~$K$ via~(2.7), and so this equation is actually nonlinear. 
Explicitly, it looks as 
\begin{equation}
i\hb \pa_x K_{x,y}(z)=\frac1{(2\pi\hb)^{2n}}
\iint \cH^\hb_x(z')K_{z',y'}(z)K_{x,y}(y')\,dm(z') dm(y').
\tag{2.25}
\end{equation}

One can ``resolve'' this nonlinearity by using the
representation of the $*$-product via operators~$L$ of the
left regular representation
(see general definitions, e.g., in~\cite{13}).
In the case of the phase space $\bR^{2n}$ it is well known that
\begin{equation}
f*g=f(L)g,
\tag{2.26}
\end{equation}
the operators where $L=(L_q,L_p)$ are Weyl-symmetrized and  
given precisely by 
$$
L=x+\frac{i\hb}{2}J\frac{\pa}{\pa x},\qquad 
x=(q,p)\in\bR^{2n},\quad 
J=\left(\begin{matrix}0&I\\-I&0\end{matrix}\right).
$$
So Eq.~(2.15) reads as a linear equation:
\begin{equation}
i\hb\pa_x K_{x,y}=\cH^\hb_x(L)K_{x,y},
\tag{2.27}
\end{equation}
where $\cH^\hb$ is the Ether Hamiltonian.

The derivation of the wave equations for the $*$-product kernel
demonstrated above is universal and could be produced over
curvilinear phase spaces.
All formulas, except (2.22) and (2.23), are general. 
Equation~(2.15) can be considered as
a {\it fundamental nonlinear equation determined by
the Ether Hamiltonian}. 
This equation generates the $*$-product kernel~$K$ 
starting from the Cauchy data~(2.16).
And Eq.~(2.27) can also be well established 
by means of an {\it a priori\/} computation 
of the left regular representation
(2.26) of the $*$-product algebra\footnote{The 
derivation
(2.7)--(2.27) in the case of the phase space $\bR^{2n}$ with
symplectic form $\omega=dp\wedge dq+\frac12 F(q)\,dp\wedge dq$
containing an additional ``magnetic'' summand 
was given in the author's talk in Manitoba University,
Winnipeg (July 2000).
The observation was that Eq.~(2.27) for the $*$-product kernel
is immediately generalized to such ``magnetic'' case and the
operators~$L$ are explicitly known in that case~\cite{52}.
For the properties of the kernel $K_{x,y}$ and the family $\bS_x$
and for their relations with symplectic transformations and connections
in the ``magnetic'' case, see \cite{29}.}. 
But first of all, one needs to know what the Ether Hamiltonian
$\cH^\hb$ is in general.

\section{Zero curvature equation}

In the case of a general phase space, i.e.,  
of a curvilinear symplectic manifold, 
the Ether Hamiltonian $\cH^\hb$ is a $1$-form
on this space (with respect to the ``time'' variable)
whose values are functions on the same space (with respect
to the ``space'' variable): 
\begin{equation}
\cH^\hb_x(z)=\sum^{2n}_{j=1}\cH^\hb_x(z)_j\,dx^j.
\tag{3.1}
\end{equation}

Equations~(2.12), (2.15), or (2.27) are actually systems of
equations, and so their generators $\wh{\cH}^\hb_j$
($j=1,\dots,2n$) should satisfy a compatibility condition. 
It can be written in the form 
$$
i\hb\pa \wh{\cH}^\hb=\frac12[\wh{\cH}^\hb \wedgeco \wh{\cH}^\hb],
$$
or
\begin{equation}
\pa \cH^\hb+\frac{i}{2\hb}[\cH^\hb \wedgeco \cH^\hb]\zs*=0.
\tag{3.2}
\end{equation}
Here the differential $\pa$ is taken only with respect to the 
``time'' variable and the brackets $[f,g]\zs*=f*g-g*f$ are 
taken only with respect to the ``space'' variable.

Condition (3.2) is the {\it zero curvature equation for 
the connection over the phase space determined by the Ether
Hamiltonian on the bundle whose fibers are quantum function
algebras over the same phase space.} 

In the classical limit $\hb=0$, 
Eq.~(3.2) reads\footnote{In~\cite{30}, 
where the construction of the work~\cite{6} was analyzed in the
classical limit $\hb=0$,  a zero curvature equation was obtained, 
which in a sense is analogous to~(3.3). 
An essential distinction from our equation~(3.3) is that,    
instead of the phase space Poisson brackets $\{\cdot,\cdot\}$,  
the commutator of formal vector fields 
along tangent fibers is used in \cite{30}.
On the quantum level, the so-called ``Weyl bundle'' used in
\cite{6,6a} consists of algebras of polynomials along tangent
fibers whose algebraic structure is just the Groenewold--Moyal
product, in contrast to (3.2), where the $*$-product is taken
over the whole phase manifold.}
\begin{equation}
\pa\cH+\frac12\{\cH\wedgeco\cH\}=0.
\tag{3.3}
\end{equation}
Here we omit the index~$\hb$ in the Ether Hamiltonian when 
considering the limit $\hb=0$.
The Poisson brackets in (3.3) correspond to the
symplectic form 
\begin{equation}
\omega(z)=\frac12\omega_{jk}(z)\,dz^k\wedge dz^j.
\tag{3.4}
\end{equation}

Note that the symplectic case under study 
is a particular case of the general Poisson situation.
In the general Poisson picture, 
the ``wave'' equations (2.12), (2.15) and the 
zero curvature equations (3.2), (3.3) were  
used in \cite{12,31} (see also \cite{13}) 
and the $1$-form 
$\cH$ was called the Cartan structure. 
In this situation the ``space'' argument belongs to the
Poisson manifold itself, while the ``time'' argument belongs
to  the dual manifold, i.e., to a finite-dimensional
pseudogroup. 
In the symplectic case this pseudogroup can be identified with
the ``space'' manifold, and we again obtain (2.12), (2.15), 
(3.2), (3.3). Thus one can use those ``Poisson'' results 
in order to study the solutions of our basic equations 
in the symplectic case.

The following geometric lemma has general Poisson
settings~\cite{32,33} (see also \cite{13}).

\begin{lemma}
Let $\cX$ be a symplectic manifold.
In a neighborhood of the zero section, the cotangent bundle
$T^*\cX$ admits a symplectic fibration $\ell$ over the base
$\cX$ such that  

{\rm(i)} $\ell$ is the identical map on the zero section 
$T^*_0\cX\approx\cX$, that is, 
$\ell(x,\eta)\big|_{\eta=0}=0${\rm;}

{\rm(ii)} the dual fibration $r$ given by the reflection in
momenta{\rm:}
$$
r(x,\eta)\od\ell(x,-\eta),\qquad \eta\in T^*_x\cX,
$$
is in involution with~$\ell$, that is 
$$
\{r^j,\ell^k\}\zs{T^*\cX}=0,\qquad \forall\, j,k.
$$
\end{lemma}

A symplectic fibration with properties (i), (ii) will be called
{\it reflective}.

Now, let us construct the solution of the zero curvature
equation (3.3). Note that in the symplectic case under study the
matrix $\pa\ell/\pa\eta$ is not
degenerate at least in a sufficiently small neighborhood of the
zero section. 
Thus the following equation is solvable:
\begin{equation}
\ell(x,\eta)=z\quad\Longrightarrow\quad\eta=\cH_x(z).
\tag{3.5}
\end{equation}
The symplecticity of $\ell$ means
$$
\{\ell^j,\ell^k\}\zs{T^*\cX}=\Psi^{jk}(\ell),\qquad 
\Psi\od\omega^{-1},
$$
and so it is easy to check that the solution of~(3.5) satisfies 
\begin{equation}
\pa_j\cH_k-\pa_k\cH_j+\{\cH_j,\cH_k\}\zs{\cX}=0.
\tag{3.3a}
\end{equation}
Here we use the notation $\pa_j=\pa/\pa x_j$ for the
derivatives, and the Poisson brackets are taken with respect to
$z\in \cX$. The system (3.3a) coincides with~(3.3).

From statements (i), (ii) of Lemma~3.1, we obtain
\begin{equation}
\cH_x(z)\big|_{z=x}=0,\qquad
\frac12D\cH_x(z)\big|_{z=x}=\omega(x).
\tag{3.6}
\end{equation}
Here we use the notation $D=\pa/\pa z$ for the derivatives.

For the second derivatives at $z=x$,
from~(3.3a) we obtain the formula
\begin{equation}
\frac12D^2\cH_x(z)\big|_{z=x}=\omega(x)\Gamma(x),
\tag{3.7}
\end{equation}
where $\Gamma$ is a symplectic connection on~$\cX$, that is,
$$
\pa\omega_{jk}-\Gamma^\ell_{js}\omega_{\ell k}
-\Gamma^\ell_{ks}\omega_{j\ell}=0,
\quad\text{or}\quad \nabla\omega=0.
$$

This is a crucial place: the symplectic connection $\Gamma$ has
appeared automatically from the symplectic fibration~$\ell$
of the secondary phase space (local symplectic
groupoid\footnote{The  
fibration~$\ell$ is the target mapping (or the left restriction
mapping) in the local symplectic groupoid over $\cX$ which is
realized as a neighborhood of the zero section of $T^* \cX$
(\cite{31,34}).})   
over the original manifold~$\cX$.

Of course, instead of (3.6), (3.7), one can write direct formulas
for~$\Psi=\omega^{-1}$ and for~$\Gamma$ via~$\ell$:
$$
\Psi(x)=2\frac{\pa\ell}{\pa\eta}(x,\eta)\bigg|_{\eta=0},\qquad
\Gamma(x)=4\omega(x)\cdot\frac{\pa^2\ell}{\pa\eta\pa\eta}(x,\eta)
\bigg|_{\eta=0}\cdot\omega(x).
$$

Let us fix a point $x\in\cX$ and a tangent vector $v\in T_x\cX$.
Denote by $\Exp_x(2v\tau)$ the trajectory of the Hamiltonian
$v\cH_x$
starting from the point~$x$. Then 
\begin{equation}
\cH_x\big(\Exp_x(-v)\big)=-\cH_x\big(\Exp_x(v)\big).
\tag{3.8}
\end{equation}

\begin{theorem}
Let $(\cX,\omega,\Gamma)$ be a symplectic manifold $\cX$ 
with symplectic form $\omega$ and symplectic connection~$\Gamma$. 
In a neighborhood of the diagonal $z\equiv x$,
there is a unique solution of the zero curvature
equation~\thetag{3.3}, 
which obeys conditions~\thetag{3.6}, \thetag{3.7}, and~\thetag{3.8}.
This solution can be obtained by~\thetag{3.5}
via a unique reflective symplectic fibration~$\ell$ over~$\cX$.
\end{theorem}

We call $\cH$ the {\it classical Ether Hamiltonian}.
The mapping 
$$
\Exp_x:\, T_x\cX\to \cX
$$
will be called the Ether exponential 
mapping and the trajectory $\Exp_x(2v\tau)$ will be called 
the {\it Ether geodesics\/} through the point~$x$ 
with velocity~$v$. 

Note that Ether geodesics correspond to some vertical curves
in the fiber $T^*_x\cX$, namely, to the curves 
\begin{equation}
\eta(\tau)\od \cH_x\big(\Exp_x(2v\tau)\big).
\tag{3.9}
\end{equation}
They are perpendicular to the velocity $v\in T_x\cX$ and obey
the equations 
\begin{equation}
\frac{d}{d\tau}\eta(\tau)=v\Omega^{[x]}(\eta(\tau)),\qquad
\eta(\tau)\bigg|_{\tau=0}=0. 
\tag{3.10}
\end{equation}
Here the tensor $\Omega^{[x]}$ is determined in a neighborhood
of zero in $T^*_x\cX$ by the following relations
\begin{equation}
\Omega^{[x]}_{lk}(\cH_x)=\{\cH_{xl},\cH_{xk}\}.
\tag{3.11}
\end{equation}
This tensor introduces a symplectic structure to the fiber
$T^*_x\cX$ so that the {\it Ether mapping $\cH_x$ is symplectic}. 
With respect to this structure, system (3.10) is Hamiltonian and
the related Hamilton function is just the linear function~$v$
(i.e., the function $\langle v,\eta\rangle$, $\eta\in T^*_x\cX$).

Besides other properties,
let us stress that the infinitesimal geometry of space,
including its symplectic structure, connection, and curvature, is
sitting inside the Ether Hamiltonian 
which can be considered as a ``generating function'' 
for this kinematic geometry.
Using the Ether exponential coordinates $ z=\Exp_x(v)$, we represent
the components of this ``generating function'' as follows:
\begin{equation}
\cH_j=2\omega_{jk}(x)v^k
+2\omega_{km}(x)R^m_{l,js}(x)v^kv^lv^s+O(v^5),
\tag{3.12}
\end{equation}
where $R$ is the curvature of the symplectic connection~$\Gamma$.
The skew-symmetry with respect to the
tangent argument~$v$ corresponds to the cyclicity
condition~(3.8).  Higher terms in~(3.9) contain derivatives of
the curvature~$R$.
\medskip

Now we demonstrate an alternative way to obtain the Ether Hamiltonian. 

Let $x\in\cX$.
A symplectic mapping $s_x:\,\cX\to\cX$ is called a {\it
reflection\/}\footnote{In 
\cite{35} such a mapping is called a symmetry under the
additional condition that the fixed point is unique.
However, the term ``symmetry'' in the classical theory of
symmetric spaces carries the strong property $s_x s_y s_x=s_{s_x(y)}$, 
which does not hold in our case.} 
in~$x$ if it is an involution: $s^2_x=\id$, 
and~$x$ is an isolated fixed point: $s_x(x)=x$.
A symplectic manifold~$\cX$ with a smooth family of reflections
$\{s_x\mid x\in\cX\}$ will be called a {\it reflective\/}
symplectic manifold. 

\begin{theorem}
Let $(\cX,\omega)$ be a simply connected reflective symplectic
manifold. Then there is a natural symplectic connection on~$\cX${\rm:}
\begin{equation}
\Gamma(z)=-\frac12 D^2 s_x(z)\bigg|_{x=z},\qquad z\in\cX.
\tag{3.13}
\end{equation}
The Ether Hamiltonian corresponding to $(\cX,\omega,\Gamma)$ is
given by 
\begin{equation}
\cH_x(z)=\int^z_x \langle \omega(z)\pa s_x
\big(s_x(z)\big),\,dz\rangle.
\tag{3.14}
\end{equation}

And vise versa{\rm:} the Ether Hamiltonian corresponding to 
$(\cX,\omega,\Gamma)$ uniquely determines the reflective
structure over $\cX$ by solving the Cauchy problem
\begin{equation}
\pa s_x+\{\cH_x,s_x\}=0,\qquad
s_x(z)\bigg|_{x=z}=z,
\tag{3.15}
\end{equation}
or just by solving the implicit equation 
\begin{equation}
\cH_x\big(s_x(z)\big)=-\cH_x(z)
\tag{3.16}
\end{equation}
{\rm(}in a domain where $D\cH$ not degenerate{\rm)}.
This family of reflections is related to the connection~$\Gamma$
via formula~\thetag{3.13}.
The Ether geodesics are symmetric with respect to reflections
\begin{equation}
s_x\big(\Exp_x(v)\big)=\Exp_x(-v)
\tag{3.17}
\end{equation}
{\rm(}but, in general, reflections $s_x$ are not affine with
respect to~$\Gamma${\rm)}.
\end{theorem}

In particular, this theorem can be applied to symplectic
symmetric manifolds \cite{36}, i.e., to the case where the
curvature of~$\Gamma$ is covariantly constant ($\nabla R=0$).  
In this case the connection (3.13) is just the canonical
Cartan--Loos connection, 
the axiom $s_x s_y s_x=s_{s_x(y)}$ holds, all reflections $s_x$
are affine, and the Ether geodesics coincide with the usual
geodesics of the connection~$\Gamma$.  

As an example, let us consider 
the sphere $\bS^2=\{|x|=1\}$ embedded in
$\bR^3$ and endowed with standard symplectic form.
In this case
\begin{equation}
\cH_x(z)=2[x\times z]\cdot dx\bigg|_{\bS^2},
\tag{3.18}
\end{equation}
where the brackets $[\,\cdot\times\cdot\,]$ stay 
for the vector-product in $\bR^3$.  

\section{Intrinsic dynamical objects}

Let us denote by $g_{x,y}$ the symplectic transformations
$$
z\to g_{x,y}(z),\qquad z\in\cX,
$$ 
obtained by shifts along the Ether dynamical system 
(that is, the ``time'' derivative $\pa_x g_{x,y}(z)$
coincides with the Ether Hamiltonian vector field at the point
$g_{x,y}(z)$, and the ``initial'' data are $g_{y,y}(z)=z$).

These transformations obey the groupoid rule
\begin{equation}
g_{x,y}\cdot g_{y,z}=g_{x,z}\,,
\tag{4.1}
\end{equation}
and produce a certain group of intrinsic transformations of the
phase space.
This group, in general, has an infinite dimension,
and we prefer to speak about the finite-dimensional 
{\it groupoid of Ether translations}. 

In terms of symplectic reflections $s_x$ 
the symplectic transformations $g_{x,y}$ are expressed by the
formula 
\begin{equation}
g_{x,y}=s_x\cdot s_y,
\tag{4.2}
\end{equation}
and so
$$
s_x(y)=g_{x,y}(y).
$$
\medskip

Other important transformations are related to the Ether
exponential mapping. Let us fix $x\in\cX$, 
and let $v\in T_x\cX$.
Recall that by $\Exp_x(2vt)$ we denote the trajectory of the
Hamiltonian $v\cH_x$ starting from the
point~$x$. But there are other trajectories starting from other
points. We consider all of them. 

Denote by $e^{2vt}_x(y)$ the trajectory of the Hamiltonian 
$v\cH_x$ starting from the point~$y$.
In particular, if $y=x$, then we have

\begin{equation}
e^v_x(x)\equiv \Exp_x(v).
\tag{4.3}
\end{equation}
\medskip

\begin{lemma}
For each $x\in\cX$ the pseudogroup of symplectic transformations
$\{e^v_x\}$ obeys the identity
\begin{equation}
s_x\cdot e^v_x=e^{-v}_x\cdot s_x.
\tag{4.4}
\end{equation}
\end{lemma}

We call $e^v_x$ {\it exponential transformations}.
\medskip

What is the quantum version of all these transformations?

Heuristically, we associate the reflections $s_x$ 
with the operators $\bS_x$, 
the translations $g_{y,x}$ with the operators 
$\bS_x\bS_y=\wh{K}_{x,y}$, 
and the exponential transformations $e^v_x$ 
with the operators $\exp\{\frac{i}{2\hb}v\wh{\cH}^\hb_x\}$.

More precisely, we introduce quantum mappings 
$\wh{s}^*_x$, $\wh{g}^*_{y,x}$, and $\wh{e}^{\,v*}_x$, 
acting in the $*$-product algebra, 
by evaluating the Heisenberg transforms. 
Say, we could define the mapping $f\to\wh{s}^*_x f$ 
in the following way:
\begin{equation}
\wh{\wh{s}^*_x f}=\bS_x\wh{f}\,\bS^{-1}_x,
\tag{4.5}
\end{equation}
or
\begin{equation}
(\wh{s}^*_x f)(z)=\tr\big(\bS^{-1}_x\bS_z\bS_x\wh{f}\big).
\tag{4.5a}
\end{equation}

If one does not know the operator realization of the
$*$-product algebra,
then, instead of operators $\bS_x$,
it is possible to use the differential equations with respect to
the parameters $x,y,\dots\,$.
These equations follow form our basic relations (2.12)
and~(2.15). 

The solution of the Cauchy problem 
\begin{equation}
\cases
\ds \pa F+\frac{i}{\hb}[\cH^\hb,F]_*=0,
\\
F\big|_{\text{diagonal}}=f,
\endcases
\tag{4.6}
\end{equation}
we denote by $F(x,z)=(\wh{s}^*_xf)(z)$, 
and we call the mappings $\wh{s}^*_x$ {\it quantum reflections}.

In Eqs.~(4.6) the unknown function $F=F(x,z)$ 
depends on $x,z\in\cX$, 
the differential~$\pa$ is taken with respect to~$x$, 
the quantum brackets $[\,\cdot,\cdot\,]_*$ 
are taken with respect to~$z$, 
and the Cauchy data are given on the diagonal $x=z$. 

If we consider the same equation (4.6) but for the function
$\cF=\cF(x,y,z)$ depending on the additional argument $y\in\cX$ 
and change the Cauchy data by the following ones:
\begin{equation}
\cF(x,y,z)\big|_{x=y}=f(z),
\tag{4.7}
\end{equation}
then the solution $\cF$ determines the {\it quantum
translations}~$\wh{g}^*_{y,x}$ as 
\begin{equation}
\cF(x,y,z)=(\wh{g}^*_{y,x}f)(z).
\tag{4.8}
\end{equation}

At last, the solution of the Cauchy problem
\begin{equation}
\frac{\pa\cE}{\pa t}+\frac{i}{\hb}[\cE, v\cH^\hb_x\,]_*=0,
\qquad 
\cE\bigg|_{t=0}=f,
\tag{4.9}
\end{equation}
taken at $t=1/2$,
determines the {\it quantum exponentials}
\begin{equation}
\cE\big|_{t=1/2}=\wh{e}^{\,v*}_x f,\qquad
v\in T_x\cX.
\tag{4.10}
\end{equation}

The restriction of this function to the diagonal (as in (4.3))
can be called a {\it quantum Ether exponential mapping\/}: 
\begin{equation}
f\to f\big|_{\wh{\Exp}_x(v)} \od (\wh{e}^{\,v*}_x f)(x).
\tag{4.11}
\end{equation}

On the quantum level we have analogs of identity (4.4):
\begin{equation}
\wh{e}^{\,v*}_x \cdot \wh{s}^{\,*}_x 
=\wh{s}^{\,*}_x \cdot \wh{e}^{\,-v*}_x,
\tag{4.12}
\end{equation}
and of identities (3.8), (3.16):
\begin{align}
\wh{s}^{\,*}_x \cH^\hb_x&=-\cH^\hb_x,
\tag{4.13}\\[10pt]
\cH^\hb_x\big|_{\wh{\Exp}_x(-v)}&=-\cH^\hb_x\big|_{\wh{\Exp}_x(v)}.
\tag{4.14}
\end{align}
\medskip

All these formulas assume that one knows the $*$-product, 
and so it is possible to solve Eqs.~(4.6) and (4.9) at least
asymptotically as $\hb\to0$.

Under the usual assumption
\begin{equation}
f*g\simeq fg+\sum_{k\geq1}\hb^k c_k(f,g),
\tag{4.15}
\end{equation}
where $c_k$ are bidifferential operators of order~$k$, 
\begin{equation}
c_1=-\frac{i}{2}\{\,\cdot,\cdot\,\},\qquad 
c_k(f,g)=(-1)^k c_k(g,f)=\ol{c}_k(g,f),
\tag{4.16}
\end{equation}
the quantum mappings defined above correspond to 
their classical counterparts:
\begin{equation}
\begin{aligned}
\wh{s}^*_x&=s^*_x\big(I+O(\hb^2)\big),\\
\wh{e}^{\,v*}_x&=e^{v*}_x\big(I+O(\hb^2)\big),\\
f\big|_{\wh{\Exp}_x(v)}&=\big(I+O(\hb^2)\big) f\big|_{\Exp_x(v)}.
\end{aligned}
\tag{4.17}
\end{equation}
Here the quantum corrections $O(\hb^2)$ are formal
$\hb^2$-series of the following type:
$$
O(\hb^2)=\hb^2a_1+\hb^4 a_2+\dots ;
$$
their coefficients $a_k$ are differential operators of order
$2k+1$ easily derived from (4.6), (4.9) via the
coefficients~$c_k$~(4.15).  

In the next section we consider the procedure of constructing the
coefficients~$c_k$.

\section{Weak asymptotics and quantum zero curvature equation}

We begin to study the basic equation (2.15) in order to
construct the $*$-product over~$\cX$ 
at least asymptotically as $\hb\to0$. 
There are two types of asymptotic expansions of $*$-products: 
weak and semiclassical. 
Both are useful and both follow from the exact
fundamental nonlinear equation (2.15) for the $*$-product kernel
given by the Ether Hamiltonian.

Note that the semiclassical view on the quantization problem
over general symplectic manifolds has been developed 
from the first papers \cite{37,38,4} 
(summarized in the book \cite{13}). 
In \cite{4} this method was called {\it asymptotic
quantization}. 
The $*$-product was constructed there by matching local
semiclassical expansions, and the quantization condition 
\begin{equation}
\frac1{2\pi\hb}[\omega]-\frac12 c_1\in H^2(\cX,\bZ)
\tag{5.1}
\end{equation}
first appeared as a necessary condition for the
existence of an irreducible operator representation of 
the $*$-product~\cite{38,4,13}.

The weak asymptotics approach to the quantization problem, 
pioneered in the fundamental papers \cite{1,39}, was a
basis of all those semiclassical developments. 
In \cite{4} the name 
{\it deformation quantization} was proposed for 
the weak asymptotics method. 
The general results in this deformation direction were obtained
in \cite{40}--\cite{43}, \cite{6,6a}. 

Now we apply our basic evolutionary equations to construct 
a unique $*$-product in the standard deformation form (4.15):
\begin{equation}
(f*g)(z)=\Big(1+\sum_{k\geq1}\hb^k c_k\Big)f(x)g(y)
\Big|_{x=y=z}.
\tag{5.2}
\end{equation}
Here $c_k$ are certain differential operators acting by the
$x,y$-arguments and obeying conditions~(4.16). 
The corresponding weak asymptotics of the integral kernel is
\begin{equation}
K_{x,y}=(2\pi\hb)^{2n}\Big(1+\sum_{k\geq1}\hb^k c'_k\Big)
\delta_x\delta_y,
\tag{5.3}
\end{equation}
where the prime $'$ means the transposition with respect to the
measure~$dm$. 

The easiest way to find the coefficients $c_k$ in~(5.2) or~(5.3)
is to substitute expression~(5.2) into~(2.20) and (4.14).
Then we uniquely derive all $c_k$ in terms of the Taylor
expansion~(3.12) 
of the Ether Hamiltonian as follows:
\begin{equation}
f*g=fg-\frac{i\hb}{2}f\underset{\leftarrow}{\nabla} \Psi
\underset{\rightarrow}{\nabla}g
-\frac{\hb^2}{8}f(\underset{\leftarrow}{\nabla}\Psi
\underset{\rightarrow}{\nabla})^2g+\ldots
\tag{5.4}
\end{equation}
Here $\Psi=\omega^{-1}$ is the Poisson tensor on the phase space, 
$\nabla$ is the covariant derivative corresponding to the
symplectic connection~$\Gamma$,
the lower arrows ($\leftarrow$ or $\rightarrow$) 
indicate the multiplier ($f$ or $g$) on which 
the derivative acts. 
Higher terms of (5.4) contain the curvature of~$\Gamma$.

Note that, 
together and simultaneously with deriving formulas (5.2), (5.4),
we have to satisfy the quantum zero curvature equation (3.2)
in which the same $*$-product (5.2) should be used. 
Thus we are looking for the quantum Ether Hamiltonian in the form
\begin{equation}
\cH^\hb=\cH+\hb^2 \cC+\hb^4\cD+\dots,
\tag{5.5}
\end{equation}
where $\cH$ is the classical Ether Hamiltonian and the quantum
corrections $\cC, \cD,\dots $ obey  
variations of the  zero curvature equation.

\begin{theorem}
Let $(\cX,\omega,\Gamma)$ be a symplectic simply-connected
manifold, and let $\omega^\hb=\omega+\hb^2\alpha+\hb^4\beta+\dots$ 
be a closed deformation of $\omega$.
Then the Ether structure generates uniquely,
explicitly, and geometrically 
{\rm(}via $\omega^\hb$ and $\Gamma${\rm)} 
the formal $\hb$-power expansions of the $*$-product 
\thetag{5.2}, \thetag{5.4} and of the solution \thetag{5.5} 
of the quantum zero curvature
equation near the diagonal.
On the diagonal the boundary condition holds:
$$
\frac12 D\cH^\hb_x(z)\bigg|_{z=x}=\omega^\hb(x).
$$
\end{theorem}

\section{Semiclassical asymptotics and quantum lift}

The weak asymptotics (5.2)--(5.4)
loses some important information about the ``wave function'' 
$K_{x,y}$. A more accurate approximation is the semiclassical
one, say, the WKB-approximation
\begin{equation}
K_{x,y}=e^{\frac{i}{\hb}\Phi_{x,y}}\varphi_{x,y}+O(\hb),
\tag{6.1}
\end{equation}
which works well outside the focal set
(in the semiclassically-simple domain). 
The difference between
(6.1) and (5.3) is the same as the difference between
(2.2) and~(2.3).

The phase $\Phi_{x,y}$ in (6.1) carries 
a dynamic geometry  of the $*$-product. 
At the ``initial moment'' $x=y$ this phase is zero, 
but its ``time'' derivatives are not zero:
\begin{equation}
\Phi_{y,y}(z)=0,\qquad 
\pa_y \Phi_{x,y}(z)\bigg|_{x=y}=\cH_y(z).
\tag{6.2}
\end{equation}
Thus, in semiclassically-simple domains,
the Ether structure is automatically generated by 
the WKB-phase of the $*$-product integral kernel. 

The $*$-product itself, after the substitution of (6.1)
into (2.7), reads
\begin{equation}
(f*g)(z)=\frac1{(2\pi\hb)^{2n}}\iint 
\Big(e^{\frac{i}{\hb}\Phi_{x,y}(z)}\varphi_{x,y}(z)+O(\hb)\Big)
f(x)g(y)\,dm(x) dm(y).
\tag{6.3}
\end{equation}

In order to compute the semiclassical approximation (6.1) 
of the $*$-product kernel, we use the Ether wave equation and
Ether translations. 

First, we resolve the nonlinearity of Eqs.~(2.15), (2.25)
following the scheme (2.26), (2.27).
For this, one needs to know operators of the left regular
representation. We compute them in the same way as in
\cite{12,13}. 

\begin{lemma}
The left regular representation of the $*$-product algebra is
given by the formula 
\begin{equation}
(f*g)(x)=f^\#(\overset{2}{x},\,-i\hb\overset{1}{\pa}_x)g(x),
\tag{6.4}
\end{equation}
where the function $f^\#$ on the secondary phase space $T^*\cX$
is determined by the equation 
\begin{equation}
f^\#(\overset{2}{x},\,-i\hb\pa_x\overset{1}{+}\cH^\hb_x*\,)1=f,
\qquad \forall x.
\tag{6.5}
\end{equation}
Here the integers $1,2,\dots$ over the operators indicate 
{\rm(}as in \cite{45}{\rm)} 
the mutual order of action of these operators.
\end{lemma}

\begin{corollary}
The asymptotic solution of Eq.~\thetag{6.5} as $\hb\to0$ is
given by 
\begin{align}
f^\#(x,\eta)
&=f(\ell(x,\eta))
-\frac{i\hb}{2}\frac{\pa^2}{\pa x^k\pa\eta_k}f(\ell(x,\eta))
\tag{6.6}\\
&\qquad 
-\frac{i\hb}{2}\frac{\pa}{\pa\eta_s}f(\ell(x,\eta))
\frac{\pa}{\pa\eta_k}\Big(\pa_k\cH_x(\ell(x,\eta))_s\Big)+O(\hb^2).
\nonumber
\end{align}
\end{corollary}

Thus the operation $\#$ in the classical limit $\hb=0$ is just
the lift of functions from $\cX$ to $T^*\cX$ by means of the
symplectic fibration~$\ell$ (see Section~3 above).
We denote $f^\#=\wh{\ell}^{\,*}f$ and call the mapping
$\wh{\ell}^{\,*}$ a {\it left quantum lift}. 
\medskip

Now let us return to our basic wave equations (2.25), (2.15)
and to apply (6.4). The result is the following.

\begin{theorem}
{\rm(i)}
The integral kernel of the $*$-product satisfies the linear
equation 
\begin{equation}
i\hb\pa_x K_{x,y}(z)
=\cL^\hb_x(\overset{2}{z},\,-i\hb\overset{1}{\pa}_z)K_{x,y}(z).
\tag{6.7}
\end{equation}
Here the $1$-form $\cL^\hb$ over $\cX$, taking values in
functions over $T^*\cX$, is defined by means of the left quantum
lift of the quantum Ether Hamiltonian{\rm:} 
$\cL^\hb=\wh{\ell}^{\,*}\cH^\hb$, i.e., 
by procedure \thetag{6.5}{\rm:}
\begin{equation}
\cL^\hb_x(\overset{2}{z},\,-i\hb\pa_z\overset{1}{+}\cH^\hb_z*\,)1
=\cH^\hb_x,\qquad \forall x,z.
\tag{6.8}
\end{equation}
The asymptotics of $\cL^\hb$ as $\hb\to0$ is 
\begin{equation}
\cL^\hb=\cL
-\frac{i\hb}{2}\frac{\pa^2\cL}{\pa z^k\pa \xi_k}
-\frac{i\hb}2 \cA+O(\hb^2),
\tag{6.9}
\end{equation}
where 
\begin{equation}
\cL_x=\ell^*\cH_x,\quad\text{or}\quad
\cL_x(z,\xi)=\cH_x(\ell(z,\xi)),\quad
\xi\in T^*_z\cX,
\tag{6.10}
\end{equation}
and
\begin{equation}
\cA_x(z,\xi)_j
=\frac{\pa}{\pa\xi_s}\cL_x(z,\xi)_j
\frac{\pa}{\pa\xi_k}\pa_k\cH_z(\ell(z,\xi))_s.
\tag{6.11}
\end{equation}

{\rm(ii)} If the arguments $x,y,z$ are considered all together 
without separation of ``space'' and ``time'' variables, 
then the linear equation for the kernel $K=K_{x,y}(z)$ reads
\begin{equation}
i\hb\,dK=\Delta^\hb K.
\tag{6.12}
\end{equation}
Here
$$
\Delta^\hb 
=\cL^\hb_x(\overset{2}{z},\,-i\hb\overset{1}{\pa}_z)
+\cL^\hb_y(\overset{2}{x},\,-i\hb\overset{1}{\pa}_x)
+\cL^\hb_z(\overset{2}{y},\,-i\hb\overset{1}{\pa}_y),
$$
where the $1$-form $\cL^\hb$ is defined in \thetag{6.8} and
\thetag{6.9}. 
\end{theorem}

So, for the integral kernel of the $*$-product
we derived the linear(!) pseudodifferential equation (6.7) 
equipped with the Cauchy data (2.16)
or the linear equation (6.12) with the Cauchy data
\begin{equation}
K\bigg|_{x=y=z}=\mu^2.
\tag{6.13}
\end{equation}

Now it is a routine exercise to compute the semiclassical
asymptotics of the solution $K$ in the semiclassically-simple
domain without focal effects. 
We substitute (6.1) into (6.7) or (6.12) and derive a
Hamilton--Jacobi equation for the phase $\Phi=\Phi_{x,y}(z)$ and
a transport equation for the amplitude
$\varphi=\varphi_{x,y}(z)$.  

\section{Intrinsic Hamilton--Jacobi and transport equations}

To the pseudodifferential equation (6.7),
one has to assign the following {\it intrinsic\/} Hamilton--Jacobi
equation (with multidimensional ``time''): 
\begin{equation}
\pa_x\Phi_{x,y}(z)+\cL_x(z,\pa_z\Phi_{x,y}(z))=0,\qquad
\Phi\big|_{x=y}=0,
\tag{7.1}
\end{equation}
where $\cL_x$ is given by (6.10) via the classical Ether Hamiltonian.

The Hamilton system (with multidimensional ``time''),
which is related to (7.1), 
is
\begin{equation}
\cases
\ds \pa_x z=D\cL_x/D\xi, &\quad z\big|_{x=y}=b\in\cX,\\
\ds \pa_x \xi=-D\cL_x/Dz, &\quad \xi\big|_{x=y}=0.
\endcases
\tag{7.2}
\end{equation}
Recall that $\xi$ is the momentum dual to $z$ so that 
\begin{equation}
\xi=\pa_z\Phi_{x,y}(z),
\tag{7.3}
\end{equation}
and the zero Cauchy data for $\xi$ in (7.2) follow from the zero
Cauchy data for $\Phi_{x,y}$ at $x=y$.

Let us denote
\begin{equation}
c=\ell(z,\xi).
\tag{7.4}
\end{equation}
Then for this equation we obtain from (7.2):
$$
\frac{\pa c^m}{\pa x^k}=\{\ell^n,\ell^m\} D_n \cH_x(c)_k.
$$
Since $\ell$ is symplectic $\{\ell^n,\ell^m\}=\Psi^{nm}(\ell)$,
then  
\begin{equation}
\frac{\pa c^m}{\pa x^k}=D_\ell \cH_x (c)_k\Psi^{\ell m}(c),\qquad
c\Big|_{x=y}=b.
\tag{7.5}
\end{equation}
So the point $c$ is the Ether translation of the point~$b$:
\begin{equation}
c=g_{x,y}(b).
\tag{7.6}
\end{equation}

Now from the first equation (7.2) it follows that
$$
\frac{\pa z^r}{\pa x^k} D_l\cH_z(c)_r=D_l\cH_x(c)_k,
$$
and from (6.16) we obtain 
$$
\frac{\pa c^m}{\pa x^k}=\frac{\pa z^r}{\pa x^k}D_l\cH_z(c)_r\Psi^{lm}(c).
$$
On the other hand, the refection equation (3.15) implies
$$
\frac{\pa s^m_z}{\pa x^k} 
=\frac{\pa z^r}{\pa x^k}\,\frac{\pa s^m_z}{\pa z^r}
=\frac{\pa z^r}{\pa x^k} D_l\cH_z(s_z)_r \Psi^{lm}(s_z).
$$
Comparing this relation with the previous one, we see that $c^m$
are equal to $s^m_z$ if their initial data at $x=y$ coincide.
Thus we prove that 
\begin{equation}
c=s_z(b).
\tag{7.7}
\end{equation}

By definition~(7.4), we also have 
\begin{equation}
\xi=\cH_z(c)=\cH_z(s_z(b))=-\cH_z(b).
\tag{7.8}
\end{equation}

From (7.3), (7.8), and (7.1) we get the derivatives of
$\Phi_{x,y}(z)$ with respect to $z$ and with respect to $x$. 
The derivative  with respect to $y$ can be computed, say, 
from the skew-symmetry condition (2.18) which implies 
$\Phi_{x,y}(z)=-\Phi_{y,x}(z)$.
So, we have 
$\pa_y\Phi_{x,y}(z)=-\pa_y\Phi_{y,x}(z)
=\cH_y\big(\ell(z,\pa_z\Phi_{y,x}(z))\big)
=\cH_y\big(\ell(z,-\pa_z\Phi_{y,x}(z))\big)
=\cH_y\big(\ell(z,\cH_z(b))\big)
=\cH_y(b)$.
Here we have again used (7.3), (7.8), and (7.1). 

Thus we have proved the following statements.

\begin{lemma}
{\rm(i)} The pair of Eqs.~\thetag{7.7}, \thetag{7.8}
determines the trajectories $z=z(x')$, $\xi=\xi(x')$ of the
Hamilton system \thetag{7.2} in $T^*\cX$ via the trajectories
$c=c(x')$ of the Ether system \thetag{7.5} in~$\cX$.

{\rm(ii)} In the semiclassically-simple domain\footnote{i.e., in
a domain where a solution of (7.10) exists and is unique.}
the phase of the $*$-product kernel is determined 
by the Ether Hamiltonian as follows{\rm:}
\begin{equation}
d(\Phi_{x,y}(z))=\cH_x(a)+\cH_y(b)+\cH_{z}(c).
\tag{7.9}
\end{equation}
Here, on the left-hand side,  
the differential is taken with respect to all variables $x,y,z$
and, on the right-hand side, the points $a,b,c$ are taken 
from the equations
\begin{equation}
c=s_z(b),\qquad b=s_y(a),\qquad a=s_x(c).
\tag{7.10}
\end{equation}
These points are related to one another via the
pseudogroup operation $\ostar$ on $\cX$ corresponding to
the ``Cartan structure'' as in \cite{31}:
\begin{equation}
z=x\,\overset{a}{\ostar}\,y,\qquad
x=z\,\overset{b}{\ostar}\,y,\qquad
y=z\,\overset{c}{\ostar}\,x.
\tag{7.11}
\end{equation}

{\rm(iii)} The Poincare--Cartan integral representation for the
solution of the Ha\-mil\-ton--Ja\-co\-bi equation~\thetag{7.1} is
\begin{equation}
\Phi_{x,y}(z)=\int^{x}_{y}\cH(a)-\int^{z}_{b}\cH(b).
\tag{7.12}
\end{equation}
Here the first integral is taken along an arbitrary path $\gamma$
from~$y$ to~$x$, and the second integral is taken along the
corresponding path $\{z(x')\mid x'\in\gamma\}$ from~$b$ to~$z$ 
{\rm(}see item~{\rm(i)} above{\rm)}.
\end{lemma}

\begin{corollary}
Let the triple $x,y,z$ belong to the semiclassically-simple
domain, and let $\Delta(x,y,z)$ be the triangle membrane in
$\cX$ bounded by the Ether geodesics through the points $x,y,z$ 
connecting the vertices $a,b,c$. 
The solution of the intrinsic Hamilton--Jacobi equation~\thetag{7.1} 
is given by the symplectic area of this triangle{\rm:}
\begin{equation}
\Phi_{x,y}(z)=\int_{\Delta(x,y,z)}\omega.
\tag {7.13}
\end{equation}
\end{corollary}

This statement specifies the result of \cite{35} by fixing the
choice of sides of the triangle in (7.13), and proves that
formula (7.13) {\it actually presents the phase of the $*$-product
kernel determined by Eqs}.~(2.15)--(2.19).
\medskip

Now let us calculate the amplitude 
$\varphi=\varphi_{x,y}(z)$ in the semiclassical
representation (6.1). 
It follows from (6.7), (6.9), and (2.16) that
this amplitude obeys the {\it intrinsic\/} transport equation 
\begin{equation}
\begin{gathered}
\bigg(\frac{\pa}{\pa x}+\frac{\pa\cL_x}{\pa\xi}\frac{\pa}{\pa z}
+\frac12\tr\bigg(\frac{\pa^2\cL_x}{\pa\xi\pa\xi}
\frac{\pa^2\Phi}{\pa z\pa z}
+\frac{\pa^2 \cL_x}{\pa z\pa\xi}\bigg)
+\frac12\cA_x\bigg)\varphi=0,\\
\varphi\Big|_{x=y}=\mu^2,
\end{gathered}
\tag{7.14}
\end{equation}
where $\cL_x$ is determined by (6.10) via the classical Ether
Hamiltonian. 
In this equation, instead of the argument~$\xi$,
one has to substitute the derivative of the phase $\Phi$
from~(7.3).

Solving (7.14) directly, we obtain 
\begin{equation}
\varphi_{x,y}(z)=2^n\mu^2\det
\Big(I-D(s_z\circ s_x\circ s_y)(b)\Big)^{-1/2}.
\tag{7.13}
\end{equation}

In particular, we have proved the conjecture \cite{35}.

\begin{corollary}
{\rm(i)} 
Let $\cX$ be a reflective simply connected symplectic manifold. 
In the semiclassically-simple domain the amplitude of the
$*$-product kernel \thetag{6.1} over $\cX$ is given by formula 
\thetag{7.13} 
in which $b$ is the fixed point\footnote{This is equivalent to
Eqs.~(7.10).} 
of the mapping $s_z\circ s_x\circ s_y$.
The equivalent formula for the amplitude, which does not use
reflections, is the following{\rm:}
$$
\varphi_{x,y}(z)=2^n\mu^2
\bigg(\frac{\det\pa_y D_z\Phi_{x,y}(z)\cdot\det\omega(b)}
{\det D \cH_y(b)\cdot \det D \cH_z(b)}\bigg)^{1/2}.
$$

{\rm(ii)}
The semiclassically multiple domain of Eq.~\thetag{2.25}
is formed by those triples $x,y,z\in\cX$ for which the mapping
$s_z\circ s_x\circ s_y$ has several isolated fixed points
{\rm(}say, on the sphere $\bS^2$ any nonfocal triple 
is of multiplicity~$2${\rm)}.
If the triple $x,y,z$ belongs to a semiclassically-multiple
domain, then in formula \thetag{6.1} 
one has to take a sum over all possible Ether geodesic triangles
$\Delta(x,y,z)$ and to multiply each summand by a suitable
exponential $\exp\{-i\frac{\pi}{2}m\}$, where $m\in \bZ$ is the
Maslov index on the graph of the groupoid multiplication.  

{\rm(iii)}
The focal set of the basic equation \thetag{2.25}, where the
asymptotics of the solution is singular as $\hb\to0$, 
consists of those triples $x,y,z\in\cX$ for which 
the symplectic mapping $s_z\circ s_x\circ s_y$ has nonisolated 
fixed points. 

{\rm(iv)}
If the second Betti number of the manifold $\cX$ is nontrivial,
then the existence of the global  semiclassical solution of the
basic equation \thetag{2.25} is ensured by the quantization
condition \thetag{5.1}.
\end{corollary}

Actually, the semiclassical approximation 
of the kernel $K_{x,y}(z)$, globally in the phase space, 
including a neighborhood of the focal set, 
can also be obtained by applying some version of the ``canonical
operator''~\cite{46} to Eq.~(6.7).

There are some interesting papers about the semiclassical
approximation of the $*$-product kernel in the special
symmetric case (the curvature $R$ is covariantly constant) 
\cite{17a,28,47}. 
For instance, in \cite{47} it was first demonstrated
that the focal set, where the kernel $K_{x,y}(z)$ is not of the
WKB-type (6.1), exists not only on such topologically nontrivial
manifolds like $\mathbb S^2$, but even on the Lobachevski plane
with its standard symplectic form.

The formal deformation quantization in the symmetric symplectic
case was first constructed in \cite{36}.
It would be interesting to check whether
the explicit $\hb$-power series for the $*$-product on symmetric
spaces obtained in \cite{36} follows as a particular case from our 
wave equations.

There is another special case, namely, K\"ahlerian manifolds. 
In \cite{48},  
the presence of the complex structure allowed 
deriving integral representations for $*$-prod\-ucts
via coherent states under some additional conditions. 
The most critical condition is that the Liouville measure should
be a reproducing measure (up to a constant multiplier; see \cite{48}).
This condition certainly holds on homogeneous manifolds. 
Schematically, this means that the integral $*$-product kernel
has the form $K=\exp\{\frac{i}{\hb}\Phi\}\cdot\const$, i.e., 
the amplitude $\varphi$ is constant, assuming that the 
integration measure $dm$ is the Liouville one. 
Of course, on general
K\"ahlerian manifolds $\varphi$ is not constant and the entire
series for the amplitude $\varphi+\hb\varphi^1+\dots$ can be
calculated explicitly \cite{49}.  
Thus the asymptotic quantization is geometrically well defined 
and the semiclassical expansion of the integral $*$-kernel is
explicitly known (without the focal set problem)
over any K\"ahlerian manifold. 
It is easy to see that 
the scheme (2.12), (2.15), (2.25)
works in this polarized case, although the cyclic property
(2.19) fails.

\section{Exterior quantum dynamics}

There is another opportunity to apply the intrinsic geometry
generated by the Ether Hamiltonian. Let us consider a
 Hamilton function $H$ on the phase space $\cX$ and
the corresponding ``exterior'' quantum equations 
\begin{equation}
i\hb\frac{d}{dt}G_t=H*G_t=G_t*H.
\tag{8.1}
\end{equation}

If the operator representation (2.5), (2.8) 
of the $*$-product algebra is known, then 
\begin{equation}
G_t(x)=\tr\bigg(\bS_x\cdot
\exp\bigg\{-\frac{it}{\hb}\wh{H}\bigg\}\bigg),\qquad G_0=1.
\tag{8.2}
\end{equation}
Here the exponential 
$\exp\big\{-\frac{it}{\hb}\wh{H}\big\}$ represents the solution
of the ``exterior'' Schr\"o\-din\-ger evolution equation.
So, the function $G_t$ is the symbol of the quantum evolution
operator.

To solve Eq.~(8.1) means that one has to take the 
exponential $\exp\{-\frac{it}{\hb}H\}$ and then replace the
classical multiplication determining this exponential 
by the quantum multiplication~$*$.
The last procedure can be made using the intrinsic Ether
dynamics following~(2.21). 

\begin{theorem}
Let $\cX$ be a reflective simply connected symplectic manifold. 

{\rm(i)} The Schr\"odinger evolution equation \thetag{8.1} is
equivalent to
\begin{equation}
i\hb\frac{d}{dt}G_t
=H^\#(\overset{2}{x},\,-i\hb\overset{1}{\pa}_x)G_t,\qquad G_0=1,
\tag{8.3}
\end{equation}
where $H^\#=\wh{\ell}^{\,*}H$ is the left quantum lift of $H$
defined by~\thetag{6.5}. 

{\rm(ii)} 
Denote by $\gamma^t\zs{H}$ the Hamilton flow generated by $H$
on~$\cX$. Then the semiclas\-sically-simple domain for problem
\thetag{8.3} as $\hb\to0$ consists of those $x\in\cX$ for which
the mapping $s_x\circ \gamma^t\zs{H}$ has a unique fixed point.

If this mapping has several {\rm(}finitely many{\rm)}
fixed points, then $x$ belongs to the semiclassically-multiple 
domain. 

The focal set of problem \thetag{8.3} consists of those~$x$ for
which the mapping $s_x\circ \gamma^t\zs{H}$ has nonisolated
fixed points. 

{\rm(iii)}
In the semiclassically-simple domain the asymptotics of $G_t$ 
is the following 
\begin{equation}
G_t=\exp\bigg\{\frac{i}{\hb}\int_{\Sigma_t}\omega\,\,
-\frac{it}{\hb}H\bigg\}\varphi_t+O(\hb).
\tag{8.4}
\end{equation}
Here $\Sigma_t(x)$ is a certain dynamic segment in $\cX$
{\rm(}a sickle-shaped membrane{\rm);} the exterior arc of this
segment is a Hamilton trajectory of $H$,
whose time-length is~$t$,
and the other arc of the segment connecting its ends 
{\rm(}spikes of the sickle{\rm)} is given by Ether geodesics
through the midpoint~$x$. The value of $H$ in
\thetag{8.4} is taken on the exterior arc of the segment. 

The amplitude $\varphi_t$ in formula \thetag{8.4} is given by 
\begin{equation}
\varphi_t(x)
=2^n\big(\det(I-D(s_x\circ\gamma^t\zs{H})(x_0))\big)^{-1/2}, 
\tag{8.5}
\end{equation}
where $x_0$ is the fixed point of the mapping
$s_x\circ\gamma^t\zs{H}$. 

{\rm(iv)}
In the semiclassically-multiple domain on $\cX$, the asymptotics
of $G_t$ is a sum of expressions like \thetag{8.4} over all
fixed points of the mapping $s_x\circ\gamma^t\zs{H}${\rm;}
each of the summands is multiplied by an exponential
$\exp\{-i\frac{\pi}{2}m\}$, 
where $m\in\bZ$ is the Maslov index on the graph of
$\gamma^t\zs{H}$.  
\end{theorem}

In the case of the Euclidean phase space $\cX=\bR^{2n}$,
formulas like (8.4) on 
membranes bounded by straight-line chords were first derived in 
\cite{50,51}.
In \cite{52} see also the case of the Euclidean space endowed
with the ``magnetic'' symplectic form $\omega$ and
with an additional ``electric'' form along space-time directions;
in this case membranes were constructed by means of special
``magnetic wings.''


\medskip

\begin{thebibliography}{99}

\bibitem{1}
F.~Bayen, M.~Flato, C.~Fronsdal, A.~Lichnerowicz, 
and D.~Sternheimer,  
\textit{Quantum mechanics as a deformation of classical
mechanics}, Lett. Math. Phys. \textbf{1} (1975/77),
521--530.

\bibitem{2}
F.~A.~Berezin,
\textit{General concept of quantization},
Comm. Math. Phys. \textbf{40} (1975), 153--174.

\bibitem{3}
A.~Connes,
\textit{Noncommutative Geometry},
Acad. Press, London, 1994.

\bibitem{4}
M.~V.~Karasev and V.~P.~Maslov,
\textit{Asymptotic and geometric quantization},
Uspekhi Mat. Nauk \textbf{39} (1984), no.~6, 115--173;
English transl., Russian Math. Surveys \textbf{39} (1984),
no.~6, 133--205.

\bibitem{39}
F.~Bayen, M.~Flato, C.~Fronsdal, A.~Lichnerowicz, and
D.~Sternheimer, 
\textit{Deformation theory and quantization},
Ann. Physics \textbf{111} (1978), 61--151.

\bibitem{5}
A.~Lichnerowicz,
\textit{Deformation of quantification},
Lecture Notes in Phys. \textbf{106} (1979), 209--219.

\bibitem{6}
B.~Fedosov,
\textit{Formal quantization},
in: \textit{Some Topics of Modern Math. and Their Appl. to
Problems of Math. Physics}, 1985, pp.~129--136.
(Russian).

\bibitem{6a}
B.~Fedosov,
\textit{A simple geometrical construction of deformation
quantization},
J. Diff. Geom. \textbf{40} (1994), 213--238.

\bibitem{7}
M.~Kontsevich,
\textit{Deformation quantization of Poisson menifolds},
Preprint q-alg/9709040, 1997.

\bibitem{8}
A.~Cattaneo and G.~Felder,
\textit{Poisson sigma models and symplectic groupoids}, 
math.sg/0003023.

\bibitem{9}
A.~Cattaneo and G.~Felder,
\textit{A path integral approach to the Kontsevich quantization
formula}, math/9902090.

\bibitem{10}
J.~Klauder,
\textit{Quantization is geometry. After all},
Ann. Physics \textbf{188} (1988), 120--141.

\bibitem{11}
B.~Kostant,
\textit{Quantization and representation theory},
London Math. Soc. Lecture Note Ser.
\textbf{34} (1979), 287--316.

\bibitem{46}
V.~Maslov,
\textit{Perturbation Theory and Asymptotic Methods},
Moscow Univ. Publ., 1965; 
French transl., Dunod, Paris, 1972.

\bibitem{12}
M.~V.~Karasev,
\textit{Quantization of nonlinear Lie-Poisson brackets in semiclassical
approximation}, 
Inst. Theor. Phys., Kiev,  Preprint N ITP-85-72P, 1985.
(Russian).

\bibitem{13}
M.~V.~Karasev and V.~P.~Maslov,
\textit{Nonlinear Poisson Brackets. Geometry and Quantization}, 
Nauka, Moscow, 1991;
English transl.,
Ser. Translations of Mathematical Monographs,
Vol.~119, Amer. Math. Soc., Providence, RI, 1993.

\bibitem{14}
M.~A.~Rieffel,
\textit{Deformation quantization of Heisenberg manifold}, 
Comm. Math. Phys. \textbf{122} (1989), 531--562.

\bibitem{15}
M.~V.~Karasev,
\textit{Advances in quantization: quantum tensors, explicit
star-products, and restriction to irreducible leaves},
Diff. Geometry and Its Appl. \textbf{9} (1998), 89--134. 

\bibitem{16}
M.~V.~Karasev and E.~M.~Novikova,
\textit{Non-Lie permutation relations, coherent states,
and quantum embedding}, 
in:
\textit{Coherent Transform, Quantization, and Poisson Geometry\/}
(M.~Karasev, ed.), Amer. Math. Soc., Providence, RI, 1998.

\bibitem{17}
J.~M.~Gracia-Bondia,
\textit{Generalized Moyal quantization on homogeneous symplectic
spaces},
in: \textit{Deformation Theory and Quantum Groups\/}
(M.~Gerstenhaber and J.~Stasheff, eds.),
Contemp. Math. \textbf{134} (1992), 93--114.

\bibitem{17a}
P.~Bieliavsky, 
\textit{Strict quantization of solvable symmetric spaces},
math.qa/ 0010004.

\bibitem{18}
H.~Weyl,
\textit{Reine Infinitesimalgeometrie}, 
Math. Z. \textbf{2} (1918), 384--411.

\bibitem{19}
E.~Scholz,
\textit{Hermann Weyl's ``Purely Infinitesimal Geometry,}''
in: 
\textit{Proc. Intern. Congress of Math., Z\"urich, 1994},
Birkh\"auser, Basel--Boston--Berlin, 1995,  pp.~1592--1603.

\bibitem{20}
H.~Weyl,
\textit{Selecta},
Birkh\"auser, Basel--Boston--Berlin, 1956,  p.~192.

\bibitem{21}
R.~Stratonovich, 
\textit{On distributions in representation space},
Zh. \'Exper. Teoret. Fiz. \textbf{31} (1956), 1012--1020;
English transl., Soviet Phys. JETP \textbf{4} (1957),
no.~6, 891--898.

\bibitem{22}
C.~Fronsdal,
\textit{Some ideas about quantization},
Rep. Math. Phys. \textbf{15} (1979), no.~1, 111--145.

\bibitem{23}
D.~Arual, 
\textit{The $*$-exponential},
in: \textit{Quantum Theories and Geometry}
(M.~Cahen and M.~Flato, eds.), Kluwer Akad, 1988, 23--51.

\bibitem{24}
A.~Grossmann and P.~Huguenin,
\textit{Group-theoretical aspects of the Wigner--Weyl isomorphism},
Helvetica Physica Acta \textbf{51} (1978), 252--261.

\bibitem{25}
F.~A.~Berezin
\textit{Quantization},
Izv. Akad. Nauk SSSR Ser. Mat. \textbf{38} (1974), no.~5;
English transl., Math. USSR Izv. \textbf{8} (1974), no.~5,
1109--1165.

\bibitem{26}
I.~Batalin and I.~Tyutin,
\textit{Quantum geometry},
Nuclear Phys. B \textbf{345} (1990), 645--650.

\bibitem{27}
A.~Weinstein, 
\textit{Classical theta-functions and quantum tori},
Publ. RIMS, Kyoto Univ. \textbf{30} (1994), 327--333.

\bibitem{28}
Zhao-Hui Qian,
\textit{Groupoids, Midpoints, and Quantization}, 
Thesis Univ. of California, Berkeley, 1997.

\bibitem{52}
M.~V.~Karasev and T.~Osborn,
\textit{Symplectic areas, quantization, and dynamics in 
electromagnetic fields},
J. Math. Phys.
\textbf{43} (2002), no.~2, 756--788
(quant-ph/0002041).

\bibitem{29}
M.~V.~Karasev and T.~Osborn,
\textit{Magnetic curvature of quantum phase space},
in: \textit{Proc. A.~Sakharov Conference, Phys. Inst. Russian
Akad. Sci., Moscow, June 2002} (to appear).

\bibitem{30}
C.~Emmrich and A.~Weinstein,
\textit{The differential geometry of Fedosov's quantization}, 
in: \textit{Lie Theory and Geometry. In Honor of B.~Kostant},
Progr. Math. \textbf{123} (1994), Birkh\"auser, New York, 217--240.

\bibitem{31}
M.~V.~Karasev,
\textit{Analogues of objects of the theory of Lie groups for
nonlinear Poisson brackets},
Izv. Akad. Nauk SSSR \textbf{50} (1986), no.~3, 508--538;
English transl., Math. USSR Izv. \textbf{28} (1987), no.~3,
497--527. 

\bibitem{32}
M.~V.~Karasev,
\textit{Maslov quantization conditions in higher cohomology and
analogs of notions developed in Lie theory for canonical fiber
bundles of symplectic manifolds}. I,~II.
Preprint MIEM, 1981.
Deposited at VINITI (March 12, 1982, No.~1092-82, 1093-82).
Abstract: in ``Ref. Zhurnal Matematika''(1982) No.~7A676, 7A677;
English transl., 
Selecta Math. Soviet \textbf{8} (1989), no.~3, 213--258.

\bibitem{33}
A.~Weinstein,
\textit{The local structure of Poisson manifolds},
J. Differential Geom. \textbf{18} (1983), 523--557.

\bibitem{34}
A.~Weinstein,
\textit{Symplectic groupoids and Poisson manifolds},
Bull. Amer. Math. Soc. \textbf{16} (1987), 101--104.

\bibitem{35}
A.~Weinstein,
\textit{Traces and Triangles in Symmetric Symplectic Spaces},
Contemp. Math. \textbf{179} (1994), 262--270.

\bibitem{36}
P.~Bieliavsky, M.~Cahen, and S.~Gutt,
\textit{Symmetric symplectic manifolds and deformation
quantization}, in:
\textit{Modern Group Theor. Methods in Phys.}
(J.~Bertrand et al., eds.), Kluwer Acad., 1995, pp.~63--73.

\bibitem{37}
M.~V.~Karasev and V.~P.~Maslov,
\textit{Algebras  with  general  permutation  relations  and
their applications}.  II,
Itogi Nauki i Tekhniki: Sovremennye Problemy Mat.
\textbf{13} (1979), 
VINITI, Moscow, 145--267;
English transl., J. Soviet Math. \textbf{15} (1981), no.~3, 
273--368.

\bibitem{38}
M.~V.~Karasev and V.~P.~Maslov,
\textit{Global asymptotic operators of regular
representation},
Dokl. Akad. Nauk SSSR \textbf{257} (1981), no.~1, 33--37;
English transl., Soviet Math. Dokl. \textbf{23} (1981), 
228--232.

\bibitem{40}
M.~De Wilde and P.~Lecomte,
\textit{Existence of star products and of formal deformations of
the Poisson Lie algebra of arbitrary symplectic manifolds}, 
Lett. Math. Phys. \textbf{7} (1983), 487--496.

\bibitem{41}
J.~Huebschmann,
\textit{On the quantization of Poisson algebras},
in: {\it Symplectic Geometry and Mathematical Physics,
Actes du colloque en l'honneur de J.-M.Souriau}
(P.~Donato et al., eds.)
Birkh\"auser, Basel--Boston, 1991, 204--233.

\bibitem{42}
H.~Omori, Y.~Maeda, and A.~Yoshioka,   
\textit{Weyl manifolds and deformation quantization},
Adv. Math. \textbf{85} (1991), no.~2, 224--255.

\bibitem{43}
H.~Omori, Y.~Maeda, and A.~Yoshioka,
\textit{Deformation quantization of Poisson algebras},
Contemp. Math. \textbf{179} (1994), 213--240.

\bibitem{44}
B.~Fedosov,
\textit{Deformation Quantization and Index Theory},
Akademie Verlag, Berlin, 1996.

\bibitem{45}
V.~Maslov,
\textit{Operator Methods}, 
Nauka, Moscow, 1973;
English transl., Mir, Moscow, 1976.


\bibitem{47}
G.~Tuynman and P.~Rios,
\textit{Weyl quantization from geometric quantization},
Preprint.

\bibitem{48}
M.~Cahen, S.~Gutt, and J.~Rawnsley,
\textit{Quantization of K\"ahler manifolds}. I,
J. Geom. Phys. \textbf{7} (1990), 45--62;
II, 
Trans. Amer. Math. Soc. \textbf{337} (1993), 73--98;
III,
Lett. Math. Phys. \textbf{30} (1994), 291--305;
IV,  
Lett. Math. Phys. \textbf{30} (1995), 159--168.

\bibitem{49}
M.~V.~Karasev,
\textit{Quantum surfaces, special functions, and the tunneling
effect}, Lett. Math. Phys. \textbf{59} (2001), 229--269.

\bibitem{50}
M.~Berry,
\textit{Semi-classical mechanics in phase space: a study of
Wigner's function}, 
Philos. Trans. Roy. Soc. London Ser.~A
\textbf{287} (1977), 237--271.

\bibitem{51}
M.~Marinov,
\textit{An alternative to the Hamilton--Jacobi approach in
classical mechanics}, 
J. Phys. A \textbf{12} (1979), no.~1, 31--47.

\end{thebibliography}
\end{document}